\documentclass[12pt,twoside,array,a4paper  Adobe times,
                fancyhdr,
		ulem ]{article}
\usepackage{latexsym,amsmath,amssymb,amsthm,array,
leqno,ulem,makeidx}
\begin{document}
\begin{center}
{\large{\bf {Extremal functions\\ for the anisotropic Sobolev inequalities.\\\ \\
Fonctions minimales\\ pour des in\'egalit\'es de Sobolev anisotropiques.}}}\\ \ \\
{\bf\small A. EL Hamidi({\footnotesize 1)}}\\
 {\footnotesize {\bf (1)}  Laboratoire de Math\'ematiques, 
Universit\'e de La Rochelle\\
Av. Michel Cr\'epeau, 17042 LA ROCHELLE Cedex 09
- France} \\
{\bf\small  J.M. Rakotoson({\footnotesize 2)}}\\
 {\footnotesize {\bf (2)}   Laboratoire de Math\'ematiques - U.M.R. 6086  - Universit\'e de Poitiers -
  SP2MI -\\ Boulevard Marie et Pierre Curie, T\'el\'eport 2\\ BP30179
  86962
FUTUROSCOPE CHASSENEUIL Cedex - France.}\\
\end{center}
\def\LEQ{\leqslant}
\def\leq{\leqslant}
\def\GEQ{\geqslant}
\def\geq{\geqslant}
\def\WH{\widehat}
\def\DST{\displaystyle}
\def\HF{\hfill$\diamondsuit$}
\def\R{\mathbb{R}}
\def\N{I\!\! N}
\def\L{I\!\! L}
\def\Q{I\!\!\!Q}
\def\calL{{\cal L}}
\def\calD{{\cal D}}
\def\calP{{\cal P}}
\def\calO{{\cal O}}
\def\calB{\hbox{$\cal B$}}
\def\calE{\hbox{$\cal E$}}
\def\calI{\hbox{$\cal I$}}
\def\calM{\hbox{${\mathfrak M}$}}
\def\calL{\hbox{$\cal L$}}
\def\calR{\hbox{$\cal R$}}
\def\calS{\hbox{$\cal S$}}
\def\ra{{R_{0}}}
\def\vp{{\overrightarrow p}}
\def\a{\alpha}
\def\b{\beta}
\def\eps{\varepsilon}
\def\f{\varphi}
\def\O{\Omega}
\def\D{\Delta}
\def\F{\Phi}
\def \PRF{{\noindent\bf Proof.\\ }}
\newcommand{\norm}[1]{\left\vert #1\right\vert}
\newcommand{\Norm}[1]{\left\Vert #1\right\Vert}
\newcommand{\crochet}[1]{\left\{ #1\right\}}
\newcommand{\mat}[1]{\left[ #1\right]}
\newcommand{\DE}[1]{\left({ #1}\right)}
\newcommand{\deriv}[2]{#1^{(#2)}}
\newcommand{\X}[1]{x_{#1}}
\newcommand{\UND}[3]{\norm{{#1}^{({#2})}}^{{#3}}}
\newcommand{\UFD}[3]{\big({#1}^{({#2})}\big)^{{#3}}}
\def\Max{\mathop{\rm{Max\;}}}
\def\Min{\mathop{\rm{Min\;}}}
\def\Inf{\mathop{\rm{Inf\;}}}
\def\infess{\mathop{\rm{inf\;ess\;}}}
\def\supess{\mathop{\rm{sup\;ess\;}}}
\def\essinf{\mathop{\rm{ess\;inf\;}}}
\def\esssup{\mathop{\rm{ess\;sup\;}}}
\def\oo{\mathop{\rm{ o(1)}}}
\def\OV{\overrightarrow}
\newtheorem{ex}{\bf Exercice }
\newtheorem{defi}{\bf Definition }
\newtheorem{prop}{\bf Proposition }
\newtheorem{theo}{\bf Theorem}
\newtheorem{proper}{\bf Property }
\newtheorem{coro}{\bf Corollary }
\newtheorem{lem}{\bf Lemma }
\newtheorem{remark}{\bf Remark}
\newtheorem{expl}{\bf Example}
{\bf Keywords :} Quasilinear problems, concentration-compactness, anisotropic Sobolev 
inequalities.\\ 
\begin{abstract}
The existence of multiple  nonnegative solutions to the anisotropic critical problem 
\begin{equation*} 
- \sum_{i=1}^{N} 
\frac{\partial }{\partial x_i} \left( \left| \frac{\partial u}{\partial x_i} \right|^{p_i-2} \frac{\partial u}{\partial x_i} \right) =  
|u|^{p^*-2} u \;\;  \mbox{in} \;\; \mathbb{R}^N
\end{equation*}
is proved in suitable anisotropic Sobolev spaces. The solutions corres\-pond to extremal functions 
of a certain best Sobolev constant. The main tool in our study is an adaptation of the well-known 
concentration-compactness lemma of P.-L. Lions to anisotropic operators. Futhermore, we show that the set of nontrival solutions $\calS$ is included in $L^\infty(\R^N)$ and is located 
outside of a ball of radius $\tau >0$ in $L^{p^*}(\R^N)$.\\
\begin{center}{\bf R\'esum\'e}\end{center}
Nous montrons l'existence d'une infinit\'e de solutions positives pour le probl\`eme anisotropique avec exposant critique. La m\'ethode consiste \`a regarder la meilleure constante d'une in\'egalit\'e du type Poincar\'e-Sobolev et \`a adapter le fameux principe de concentration-compacit\'e de P.L. Lions. De plus, on montre que l'ensemble des solutions $\calS$ est contenu dans $L^\infty (\R^N)$ et est localis\'e en dehors d'une boule de rayon $\tau>0$ dans $L^{p^*}(\R^N)$.\\
\end{abstract}
\section{Introduction.} 
In this paper, the existence of nontrivial nonnegative solutions to the anisotropic critical problem 
\begin{equation} \label{aniso}
- \sum_{i=1}^{N} 
\frac{\partial }{\partial x_i} \left( \left| \frac{\partial u}{\partial x_i} \right|^{p_i-2} \frac{\partial u}{\partial x_i} \right) =  
|u|^{p^*-2} u \;\;  \mbox{in} \;\; \mathbb{R}^N
\end{equation}
is studied, where the exponents $p_i$ and $p^*$ satisfy the following conditions
$$
p_i>1, \quad \quad  \sum_{i=1}^{N} \frac{1}{p_i} > 1,
$$
and the critical exponent $p^*$ is defined by
$$
p^*:=\frac{N}{\sum_{i=1}^{N} \frac{1}{p_i}-1}.   
$$

In the best of our knowledge, anisotropic equations with different orders of derivation in different directions, 
involving critical exponents were never studied before. 
In the subcritical case, we can refer the reader to the recent paper by I. Fragala {\it et al} \cite{FGK}.  

In the special case $p_i =2$, $i \in \{1, \ 2, \ \cdots , \ N \}$, Problem (\ref{aniso}) is reduced to the limiting equation  
arising in the famous Yamabe problem \cite{YAM}:
\begin{equation} \label{yam_lim}
-\Delta u =  u^{2^*-1},  \;\;\;  u > 0 \;\; \mbox{in} \;\; \mathbb{R}^N.
\end{equation}
Indeed, let $(M,g)$ be a $N$-dimensional Riemannian manifold and $S_g$ be the scalar curvature of the metric $g$. Consider 
a conformal metric $\widetilde g$ on $M$ defined by $\widetilde g:=u^{\frac{4}{N-2}} g$ whose scalar curvature 
(which is assumed to be constant) is denoted by $S_{\widetilde g}$, where 
$u$ is a positive function in $C^{\infty}(M,\mathbb{R})$. The unknown function $u$ satisfies then   
\begin{equation} \label{yamabe}
-\Delta_g u + \frac{N-2}{4(N-1)} S_g u = \frac{N-2}{4(N-1)} S_{\widetilde g} u^{2^*-1}, \;\;\; u>0
\;\;  \mbox{in} \;\; M, 
\end{equation}
where $\Delta_g$ denotes the Laplace-Beltrami operator. It is clear that, up to a scaling, the limiting problem of (\ref{yamabe}) 
(Equation (\ref{yamabe}) without the subcritical term $\frac{N-2}{4(N-1)} S_g u$) is 
exactly (\ref{yam_lim}). The question of existence of minimizing solutions to (\ref{yam_lim}) was completely solved 
by Aubin \cite{AUB} and G. Talenti \cite{TAL}. Their proofs are based on symmetrisation theory. Notice that this theory is not relevent 
in our context since the radial symmetry of solutions can not hold true because of the anisotropy of the operator. 

In \cite{PLL}, P.-L. Lions introduced the famous concentration-compactness lemma which constitutes a powerful tool for the study of 
critical nonlinear elliptic equations. The concentration-compactness lemma allows an elegant and simple proof of the existence 
of solutions to (\ref{yam_lim}) by minimization arguments. In the present work, we will adapt the concentration-compactness lemma 
to the anisotropic case and show that the infimum 
$$
\DST\Inf_{\norm u_{L^{p^*}\left(\mathbb{R}^N\right)}=1}
\left\{\sum_{i=1}^N \frac {1}{p_{i}}\Norm{\frac{\partial u}{\partial x_{i}}}^{p_{i}}_{p_{i}}
\right\}
$$ 
is achieved, of course, the functional space has to be specified. 

The motivation of the present work is to give a new result which can provide extremal functions associated to the critical level 
corresponding to anisotropic problems involving critical exponents. Notice that the genuine extremal functions are obtained 
by minimization on the Nehari manifold associated to the problem and the critical level is nothing than the energy of these extremal functions.

The natural functional framework of Problem (\ref{aniso}) is the anisotropic Sobolev spaces theory developed by \cite{NIK,VEN,RAK1,RAK2,Troisi}.
Then, let $\mathcal{D}^{1,\vec p}(\mathbb{R}^N)$ be the completion of the space $\mathcal{D}( \mathbb{R}^N )$ 
with respect to the norm 
$$
\Norm u_{1,\vp}  := \sum_{i=1}^{N} \Norm{ \frac{\partial u}{\partial x_i}}_{p_i}.
$$ 
It is well known that $\left(\mathcal{D}^{1,\vec p}(\mathbb{R}^N), \Norm \cdot_{1,\vp} \right)$ is a reflexive Banach space which is 
continuously embedded in $L^{p^*}\left(\mathbb{R}^N\right)$. \\
In what follows, we will assume that 
$$
p_{+}=\max\{p_1,\ p_2, \ ..., \ p_N \} < p^*,
$$
then $p^*$ is the critical exponent associated to the operator: 
$$\DST\sum_{i=1}^{N} 
\frac{\partial}{\partial x_i} \left(\left| \frac{\partial}{\partial x_i} \right|^{p_i-2} \frac{\partial}{\partial x_i} \right).
$$
The space $\mathcal{D}^{1, \vec p}(\mathbb{R}^N)$ can also be seen as 
$$
\mathcal{D}^{1, \vec p}(\mathbb{R}^N) = \left\{ u \in L^{p^*}(\mathbb{R}^N) \; : \; 
\left| \frac{\partial u}{\partial x_i} \right|  \in L^{p_i}(\mathbb{R}^N) \right\}.
$$

In the sequel, we will set $p_- = \min\{p_1,\ p_2, \ ..., \ p_N \}$,  
$p_+ = \max\{p_1,\ p_2, \ ..., \ p_N \}$ and $\overrightarrow p = (p_1, \ p_2, \ \cdots , \ p_n)$. 
Also, the integral symbol $\DST\int$ will denote $ \DST\int_{\mathbb{R}^N}$ and $\Norm \cdot_{p_i}$ will denote the usual Lebesgue norm 
in $L^{p_i}(\mathbb{R}^N)$. We denote by $\calM (\mathbb{R}^N)$ (resp. $\calM^+ (\mathbb{R}^N)$) the space of finite measures (resp. positive 
finite measures) on $\mathbb{R}^N$, and by $\Norm\cdot$ its usual norm. 

\section{Existence of extremal functions for a Sobolev type inequality}
In this paragraph, we shall prove that a certain best Sobolev constant is achieved.
\begin{theo}\label{t1}
Under the above assumptions on $p_i	,\ i=1,\ldots, N,\ N\GEQ2$, there exists at least one function 
$u\in \calD^{1,\overrightarrow p}(\R^N),\ u\GEQ0,\ u\neq 0$ :
$$-\sum_{i=1}^N \frac \partial{\partial x_{i}}
\left(\norm{\frac{\partial u}{\partial x_{i}}}^{p_{i}-2}\frac{\partial u}{\partial x_{i}}\right)
= u^{p^*-1} \;\;  \rm{in} \;\;  \calD'(\R^N).$$
\end{theo}
The proof will need two fundamental lemmas, the first one is a result due to M. Troisi \cite{Troisi}:
\begin{lem}\label{l1}{\rm{\bf (Troisi \cite{Troisi})}} \ \\
There is a constant $T_0>0$ depending only on $\overrightarrow p$ and $N$ such that :
$$T_0\Norm u_{p^*}\LEQ\prod_{j=1}^N\Norm{\frac{\partial u}{\partial x_{i}}}^{\frac 1N}_{p_{i}}\hbox{ and }
\Norm u_{p^*}\LEQ\frac 1{NT_0}\sum_{i=1}^N\Norm{\frac{\partial u}{\partial x_{i}}}_{p_{i}},$$
for all $u\in\calD^{1,\vp}(\R^N)$.
\end{lem}
The second lemma is a rescaling type result ensuring the conservation of suitable norms:
\begin{lem}\label{l2}\ \\
Let $\DST  \alpha_{i}= \frac {p^*}{p_{i}}-1,\ i=1,\ldots,N$.
For every $y\in \R^N$, $u\in\calD^{1,\overrightarrow p}(\R^N)$, 
and $\lambda > 0$, if we write $x=(x_{1},\ldots,x_{N}),\ y=(y_{1},\ldots,y_{N})$,
$v(x)\dot=u^{\lambda,y}(x)=\lambda u(\lambda^{\alpha_{1}}x_{1}+y_{1},\ldots,\lambda^{\alpha_{N}}x_{N}+y_{N})$, \\
we get
$$\Norm{u}_{p^*}=\Norm{v}_{p^*},$$
$$\Norm{\frac{\partial u}{\partial x_{i}}}_{p_{i}}
=\Norm{\frac{\partial v}{\partial x_{i}}}_{p_{i}},\ for\ i=1,\ldots,N,$$
thus, $\Norm u_{1,\vp}=\Norm{u^{\lambda,y}}_{1,\vp}$.
\end{lem}
\PRF
Noticing that $\DST\sum_{i=1}^N \alpha_{i}= p^*$, a straightforward computation with adequate changes of variables gives the result.
\begin{lem}\label{l3}\ \\
Let $S=\DST\Inf_{u\in \calD^{1,\vp}(\R^N),\ \Norm u_{p^*}=1}
\left\{\sum_{i=1}^N \frac 1{p_{i}}\Norm{\frac{\partial u}{\partial x_{i}}}^{p_{i}}_{p_{i}}
\right\}.$ Then $S>0.$
\end{lem}
\PRF
From Lemma \ref{l1}, we obtain that if $\Norm u_{p^*}=1$, then 
\begin{equation}\label{eq1}
\sum_{i=1}^N\Norm{\frac{\partial u}{\partial x_{i}}}_{p_{i}}\GEQ NT_{0}>0.
\end{equation}
Using standard argument, the infimum
$$
\DST\Inf\left\{\sum_{i=1}^N\frac1{p_{i}}a_{i}^{p_{i}},\ (a_{1},\ldots,a_{n})\in\R^N,\ \sum_{i=1}^N a_{i}\GEQ NT_{0}, \ 
a_{i}\GEQ0\right\}\dot= S_{1}
$$ 
is achieved and thus this minimum is positive. By relation (\ref{eq1}), one concludes that $S\GEQ S_{1}>0$.
\HF
\begin{coro}\label{c1l3}{\rm{\bf {of Lemma \ref{l3} (Sobolev type inequality)}}}\ \\
Let $p_{-}=\min(p_{1},\ldots,p_{N})$, $p_{+}=\max(p_{1},\ldots,p_{N})$ and $F$ be the real valued function defined by $F(\sigma)=\begin{cases}
\sigma^{p_{+}}&if\ \sigma\LEQ1,\\ \sigma^{p_{-}}&if\ \sigma\GEQ1.\end{cases}$\\
Then for every $u\in\calD^{1,\vp}(\R^N)$, one has
$$SF\big(\Norm u_{p^*}\big)\LEQ
\sum_{i=1}^N \frac 1{p_{i}}\Norm{\frac{\partial u}{\partial x_{i}}}^{p_{i}}_{p_{i}} \dot=P(\nabla u).$$
\end{coro}
\PRF
Let $u$ be in $\calD^{1,\vp}(\R^N)$. If $u=0$ the inequality is true. If $u\neq 0$, set $w=\DST\frac u{\Norm u_{p^*}}$, 
then from the definition of $S$ one has :
\begin{equation}\label{eq2}
\sum_{i=1}^N\frac1{p_{i}}\Norm{\frac{\partial w}{\partial x_{i}}}_{p_{i}}^{p_{i}}\GEQ S.
\end{equation}  
Since $t^{p_{i}}\LEQ t^{p_{+}}$ if $t>1$   and $t^{p_{i}}\LEQ t^{p_{-}}$ otherwise, the result follows 
from relation (\ref{eq2}) and the definition of $F$. \HF
\begin{remark}
Along this paragraph, we only need the inequality :
$$S\Norm u^{p_{+}}_{p^*}\LEQ P(\nabla u)\hbox{  whenever }\Norm u_{p^*}\LEQ1.$$
\end{remark}

We shall call $(\calP)$ the minimization problem
$$(\calP)\qquad\qquad \Inf_{\Norm u_{p^*}=1}
\left\{\sum_{i+1}^N\frac 1{p_{i}}\Norm{\frac{\partial u}{\partial x_{i}}}^{p_{i}}_{p_{i}}\right\}
=\Inf_{\Norm u_{p^*}=1}
\left\{P(\nabla u)\right\}.$$

Let $(u_{n})\subset\calD^{1,\vp}(\R^N)$ be a minimizing sequence for the problem ($\calP$). As in \cite{PLL} and Willem \cite{W}, we define the 
Levy concentration function:
$$Q_{n}(\lambda)=\sup_{y\in\R^N}
\int_{E(y,\lambda^{\alpha_{1}},\ldots,\lambda^{\alpha_{N}})}\norm{u_{n}}^{p^*}dx,\;\;\; \lambda > 0.$$

Here $E(y,\lambda^{\alpha_{1}},\ldots,\lambda^{\alpha_{N}})$ is the ellipse defined by 
$$
\DST\left\{z=(z_{1},\ldots,z_{N})\in\R^N,\ \sum_{i=1}^N\frac{(z_{i}-y_{i})^2}{\lambda^{2\alpha_{i}}}\LEQ1\right\}
$$ 
with $y=(y_{1},\ldots,y_{N})$ and $\alpha_{i}>0$ as in Lemma \ref{l2}. Since for every $n$, $\DST\lim_{\lambda\to0}Q_{n}(\lambda)=0$ and 
$Q_{n}(\lambda)\xrightarrow[\lambda\to+\infty]{}1$. There exists $\lambda_{n}>0$ such that $Q_{n}(\lambda_{n})=\DST\frac12$. Moreover there exists $y_{n}\in\R^N$ such that $$\DST\int_{E(y_{n},\lambda_{n}^{\alpha_{1}},\ldots,\lambda_{n}^{\alpha_{N}})}\norm{u_{n}}^{p^*}dx=\frac12.$$
Thus by a change of variables one  has for $v_n\dot= u_{n}^{\lambda_{n},y_{n}}$ :
$$\int_{B(0,1)}\norm{v_{n}}^{p^*}dx
=\frac12=\sup_{y\in\R^N}\int_{B(y,1)}\norm{v_{n}}^{p^*}dx.$$
Since $\Norm{v_{n}}_{p^*}=\Norm{u_{n}}_{p^*},\
 \DST\Norm{\frac{\partial v_{n}}{\partial x_{i}}}_{p_{i}}
 =\Norm{\frac{\partial u_{n}}{\partial x_{i}}}_{p_{i}},\ 
 P(\nabla u_{n})=P(\nabla v_{n})$ we deduce that $(v_{n})$ is bounded in $\calD^{1,\vp}(\R^N)$ and is also a minimizing sequence for $(\calP)$. We may then assume that :
\begin{itemize}
\item $v_{n}\rightharpoonup v$ in $\calD^{1,\vp}(\R^N)$,
\item $\DST\norm{\frac\partial {\partial x_{i}}(v_{n}-v)}^{p_{i}}
\rightharpoonup \mu_{i}$ in $ \calM^+(\R^N)$,
\item $\norm{v_{n}-v}^{p^*} \rightharpoonup\nu$ in $\calM^+(\R^N)$,
\item $v_{n}\to v$ a.e in $\R^N$.
\end{itemize}
We define :
\begin{eqnarray}
\mu &=& \sum_{i=1}^N\frac1{p_{i}} \mu_{i},\nonumber \\
\mu_{\infty} &=& \lim_{R\to+\infty}\overline{\lim_{n}}
\sum_{i=1}^N\frac1{p_{i}}\int_{\norm x>R}\norm{\frac{\partial v_{n}}{\partial x_{i}}}^{p_{i}}dx,\label{N.1}\\
\nu_{\infty} &=& \lim_{R\to+\infty}\overline{\lim_{n}}\int_{\norm x>R}\norm{v_{n}}^{p^*}dx. \qquad \qquad \label{N.2}
\end{eqnarray}
We start with some general lemmas. First by the Brezis-Lieb's Lemma~\cite{BL}, 
direct computations give  
the following
\begin{lem}\label{l4}
$$\norm{v_{n}}^{p^*}
\rightharpoonup \norm{v}^{p^*}
+\nu\ in\ \calM^+(\R^N).$$
\end{lem}
The lemma which follows gives some reverse H\"older type inequalities connecting the measures 
$\nu$, $\mu$ and $\mu_i$, $1 \leq i \leq N$.
\begin{lem}\label{l5}\ \\
Under the above statement, one has for all $\f\in C^\infty_{c}(\R^N)$
\begin{eqnarray*}
\left(\int\norm\f^{p^*}d\nu\right)^{\frac1{p^*}}
& \LEQ & \frac1{T_{0}}\prod_{i=1}^N
\left(\int\norm\f^{p_{i}}d\mu_{i}\right)^{\frac1{Np_{i}}}, \\
\left(\int\norm\f^{p^*}d\nu\right)^{\frac1{p^*}}
& \LEQ & p_{+}^{\frac1N+\frac1{p^*}}
\Norm \mu^{\frac1N+\frac1{p^*}-\frac1{p_{+}}}
\cdot\frac1{T_{0}}
\left(\int\norm\f^{p_{+}}d\mu\right)^{\frac1{p_{+}}}.
\end{eqnarray*}
\end{lem}
\PRF
Let $\f\in C_{c}^\infty(\R^N)$ and set $w_{n}=v_{n}-v$. 
Since $\DST\int\norm{\f_{x_{i}}}^{p_{i}}\norm{w_{n}}^{p_{i}}dx\xrightarrow[n\to+\infty]{}0$,  we then have :
\begin{equation}\label{eq3}
\lim_{n}\int\norm{\frac\partial{\partial x_{i}}(\f w_{n})}^{p_{i}}dx
=\lim_{n}\int\norm\f^{p_{i}}
\norm{\frac{\partial w_{n}}{\partial x_{i}}}^{p_{i}}dx
=\int\norm \f^{p_{i}}d\mu_{i}.
\end{equation}
Thus from  Lemma \ref{l1}, it follows that
\begin{equation} \label{eq4}
\left(\int\norm\f^{p^*}d\nu \right)^{\frac1{p^*}}
=\lim_{n}\left(\int\norm{\f w_{n}}^{p^*}dx\right)^{\frac1{p^*}}
\LEQ\frac1{T_{0}}\prod_{i=1}^N\left(\int\norm \f^{p_{i}}d\mu_{i}\right)^{\frac1{Np_{i}}}.
\end{equation} 
On the other hand, since
\begin{equation}\label{eq5}
\int\norm\f^{p_{i}}d\mu_{i}\LEQ p_{+}\int\norm\f^{p_{i}}d\mu  
\LEQ p_{+}\Norm \mu^{1-\frac{p_{i}}{p_{+}}}
\left(\int\norm\f^{p_{+}}d\mu\right)^{\frac{p_{i}}{p_{+}}}
\end{equation} 
applying the estimates (\ref{eq4}) and (\ref{eq5}) and knowing that $\DST\sum_{i=1}^N\frac1{p_{i}}=
1+\frac N{p^*}$, we deduce
$$\left(\int\norm\f^{p^*}d\nu\right)^{\frac1{p^*}}
\LEQ p_{+}^{\frac1N+\frac1{p^*}}
\Norm \mu^{\frac1N+\frac1{p^*}-\frac1{p_{+}}}
\cdot\frac1{T_{0}}
\left(\int\norm\f^{p_{+}}d\mu\right)^{\frac1{p_{+}}}.$$ 
This ends the proof. \ \HF 

We then have $\Norm v_{p^*}\LEQ1$. So if $\Norm v_{p^*}=1$ then $ v$ is an extremal 
function since $\DST P(\nabla v)\LEQ\liminf_{n} P(\nabla v_{n})=S$ and $S\LEQ P(\nabla v)$. 
Thus, we want to show that fact, by proving that if it is not true then we have a concentration of $\nu$ at a single point and therefore $v=0$.\\
\ \\
{\bf Main Lemma}$$\Norm v_{p^*}=1.$$

The remainder of this section is devoted to the proof of the main Lemma
\begin{lem}\label{l6} \ \\
If $  v\neq0$ then $$\DST\lim_{n}\Norm {v_{n}-v}^{p^*}_{p^*}=1-\Norm v_{p^*}^{p^*}<1.$$
\end{lem}
\PRF
From Brezis-Lieb's Lemma we have :
$$\lim_{n}\left(\Norm{ v_{n}}^{p^*}_{p^*}-\Norm{v_{n}-v}^{p^*}_{p^*}\right)=\Norm v^{p^*}_{p^*},$$
Since $\Norm{v_{n}}_{p^*}=1$, we derive the result.\HF
\begin{lem}\label{l7}   \ 
$$S\Norm \nu^{\frac{p_{+}}{p^*}}\LEQ\Norm\mu.$$
\end{lem}
\PRF
For large $n$, according to Lemma  \ref{l6}, we have :
$$\int\norm{v_{n}-v}^{p^*}dx\LEQ1.$$
Thus for all $\f\in C_{c}^\infty(\R^N),\ \norm \f_{\infty}\LEQ1$, it holds:
$$S\left(\int\norm \f^{p^*}\norm{v_{n}-v}^{p^*}\right)^{\frac{p_{+}}{p^*}}
\LEQ\sum_{i=1}^N\frac1{p_{i}}\int\norm\f^{p_{i}}\norm{\frac{\partial (v_{n}-v)}{\partial x_{i}}}^{p_{i}}dx+
\mbox{o}_n (1).$$
Letting $n\to+\infty$, one gets :
\begin{equation}\label{eq1000}
S\left(\int\norm \f^{p^*}d\nu\right)^{\frac{p_{+}}{p^*}}
\LEQ\sum_{i=1}^N\frac1{p_{i}}\int\norm\f^{p_{i}}d\mu_{i}\LEQ\Norm \mu.
\end{equation}
Using the density of $C_{c}^\infty(\R^N)$ in $C_{c}(\R^N)$, we get then 
$$S\left(\sup_{\f\in C_{c} (\R^N),\ \norm\f_{\infty}=1}\int\norm \f^{p^*}d\nu\right)^{\frac{p_{+}}{p^*}}\LEQ\Norm\mu,$$
that is the desired result.\HF
\begin{lem}\label{l8}
Let $\psi_{R}$ be in $C^1(\R),$ $0\LEQ\psi_{R}\LEQ1,$ $\psi_{R}=1\ if \norm x>R+1,$   $\psi_{R}(x)=0\ if \ \norm x<R$. Then for any $\gamma_{i}>0, \ i=0,\ldots,N$, the two equalities
\begin{eqnarray*}
\nu_{\infty}&=&\lim_{R\to+\infty}\overline{\lim_{n}}
\int\norm{v_{n}}^{p^*}\psi_{R}^{\gamma_{0}}dx, \\
\mu_{\infty}&=&\lim_{R\to+\infty}\overline{\lim_{n}}\sum_{i=1}^N
\frac1{p_{i}}\int\norm{\frac{\partial v_{n}}{\partial x_{i}}}^{p_{i}}
\psi_{R}^{\gamma_{i}}dx.
\end{eqnarray*}
hold true, where $\nu_{\infty}$ and $\mu_{\infty}$ are defined by (\ref{N.1}), (\ref{N.2}).
\end{lem}
\PRF
As in Willem  \cite{W}, one has :
$$\int_{\norm x>R+1}\norm{v_{n}}^{p^*}dx
\LEQ\int\norm{v_{n}}^{p^*}\psi_{R}^{\gamma_{0}} dx
\LEQ\int_{\norm x>R }\norm{v_{n}}
^{p^*}dx,$$
$$\int_{\norm x>R+1}\norm{\frac{\partial v_{n}}{\partial x_{i}}}^{p_{i}}dx
\LEQ\int\norm{\frac{\partial v_{n}}{\partial x_{i}}}\psi_{R}^{\gamma_{i}}
\LEQ\int_{\norm x>R}\norm{\frac{\partial v_{n}}{\partial x_{i}}}^{p_{i}}dx.$$
We conclude with the definition of $\nu_{\infty}$ and $\mu_{\infty}$.\HF
\begin{lem}\label{l9}\ \\
Let $w_{n}=v_{n}-v$. 
Then, for any $\gamma_{i}>0,\ i=0,\ldots,N$, we get
$$\nu_{\infty}
=\lim_{R\to\infty} \overline{\lim_{n} }\int\norm{w_{n}}^{p^*}
\psi_{R}^{\gamma_{0}}dx,$$
and
$$\mu_{\infty}=\lim_{R\to\infty}\overline{\lim_{n}}
\int\norm{\frac{\partial w_{n}}{\partial x_{i}} }^{p_{^i}}
\psi_{R}^{\gamma_{i}}dx.$$
\end{lem}
\PRF
Since $$\lim_{R\to+\infty}\int\norm v^{p^*}\psi_{R}^{\gamma_{0}}
=\lim_{R\to+\infty}\int\norm{\frac{\partial v}{\partial x_{i}}}^{p_{i}}\psi_{R}^{\gamma_{i}}dx=0.$$
Thus
$$\lim_{R\to\infty}\overline{\lim_{n}}\int\norm{w_{n}}^{p^*}\psi_{R}^{\gamma_{0}}dx
=\lim_{R\to\infty}\overline{\lim_{n}}\int\norm{v_{n}}^{p^*}\psi_{R}^{\gamma_{0}}dx
=\nu_{\infty}$$
and
$$\lim_{R\to\infty}\overline{\lim_{n}}\sum_{i=1}^N\frac1{p_{i}}\int
\norm{\frac{\partial w_{n}}{\partial x_{i}}}^{p_{i}}\psi_{R}^{\gamma_{i}}dx
=\lim_{R\to\infty}\overline{\lim_{n}}\sum_{i=1}^N\frac1{p_{i}}\int
\norm{\frac{\partial v_{n}} {\partial x_{i}}}^{p_{i}}\psi_{R}^{\gamma_{i}}dx.$$
\ \HF
\begin{lem}\label{l10}\ \\
$$S\nu_{\infty}^{\frac{p_{+}}{p^*}}\LEQ\mu_{\infty}.$$
\end{lem}
\PRF
From Lemma \ref{l6}, we know that for $n$ large enough, we have
$$\int\psi_{R}^{p^*}\norm{w_{n}}^{p^*}\LEQ\int\norm{w_{n}}^{p^*}dx\LEQ1.$$
Thus by Sobolev inequality (Corollary \ref{c1l3} of Lemma \ref{l3}), it follows
$$S\left(\int\norm{\psi_{R} w_{n}}^{p^*}dx\right)   ^{\frac{p_{+}}{p^*}}
\LEQ \sum_{i=1}^N\frac1{p_{i}}
\int\norm{  \frac\partial{\partial x_{i}}(\psi_{R}w_{n}) } ^{p_{i}},$$
\begin{equation}\label{eq6}
S\left(\lim_{R\to+\infty}\overline{\lim_{n}}\int\norm{\psi_{R} w_{n}}^{p^*}dx\right)   ^{\frac{p_{+}}{p^*}}
\LEQ\lim_{R\to+\infty}\overline{\lim_{n}} \sum_{i=1}^N\frac1{p_{i}}
\int\norm{  \frac\partial{\partial x_{i}}(\psi_{R}w_{n}) } ^{p_{i}}.
\end{equation}
Since
$${\lim_{n}} \sum_{i=1}^N\frac1{p_{i}}
\int\norm{  \frac{\partial\psi_{R}}{\partial x_{i}}}^{p_{i}}\norm{w_{n} } ^{p_{i}}=0,$$
then
$$\lim_{R\to+\infty}\overline{\lim_{n}} \sum_{i=1}^N\frac1{p_{i}}
\int\norm{  \frac\partial{\partial x_{i}}(\psi_{R}w_{n})}^{p_{i}}
=\lim_{R\to+\infty}\overline{\lim_{n}} \sum_{i=1}^N\frac1{p_{i}}
\int\norm{  \frac{\partial w_{n}}{\partial x_{i}}}^{p_{i}}\psi_{R} ^{p_{i}}=\mu_{\infty}.$$
relation (\ref{eq6}) and Lemma \ref{l9} give :
$$S\nu^{\frac{p_{+}}{p^*}}_{\infty}\LEQ\mu_{\infty}.$$ \ \HF

Following again the arguments used in \cite{W} we claim that:
\begin{lem}\label{l11}
$$1=\lim_{n}\Norm{v_{n}}^{p^*}_{p^*}=\Norm v^{p^*}_{p^*}+\Norm \nu+\nu_{\infty}.$$
\end{lem}
\PRF
From Lemma \ref{l4}, we have :
$$\norm{v_{n}}^{p^*}
\rightharpoonup\norm v^{p^*}
+\nu.$$
Thus 
$$\lim_{R\to+\infty}{\lim_{n}}\int(1-\psi_{R}^{p^*})\norm{v_{n}}^{p^*}dx
=\int\norm v ^{p^*}dx+\int d\nu.$$
Rewriting $\Norm {v_{n}}^{p^*}_{p^*}$ as
$$\Norm {v_{n}}^{p^*}_{p^*}
=\int(1-\psi_{R}^{p^*})\norm{v_{n}}^{p^*}
+\int\psi_{R}^{p^*}\norm{v_{n}}^{p^*},
$$
we obtain
\begin{eqnarray*}
\lim_{n}\Norm {v_{n}}^{p^*}_{p^*}
&=&\lim_{R\to+\infty}\lim_{n}\int(1-\psi_{R}^{p^*})\norm{v_{n}}^{p^*}
+\lim_{R\to+\infty}\overline{\lim_{n}}\int\psi_{R}^{p^*}\norm{v_{n}}^{p^*}\\
&=&\Norm v^{p^*}_{p^*}+\Norm\nu+\nu_{\infty}
\end{eqnarray*}
\ \HF\\
Next, we shall prove the following corollary:
\setcounter{coro}{0}
\begin{coro}\label{c1l5}{\rm{\bf{(of Lemma \ref{l5})}}}\ \\
There exists an at most countable index set $J$ of
 distinct  points $\{ x_{j}\}_{j\in J}\subset\R^N$ and nonnegative weights 
$a_{j} \hbox{ and }b_{j},\ j\in J$  such that :
\begin{enumerate}
\item $\DST\nu=\sum_{j\in J}a_{j}\delta_{x_{j}}$.
\item $\DST\mu\GEQ\sum_{j\in J}b_{j}\delta_{x_{j}}$.
\item 
 $Sa_{j}^{\frac{p_{+}}{p^*}}\LEQ b_{j},\ \forall j\in J$.
\end{enumerate}
\end{coro}
\PRF
The proof follows essentially the concentration compactness principle of P.L. Lions \cite{PLL} because we have the reverse H\"older type inequalities of Lemma \ref{l5}.
 
Indeed, the second statement of this lemma implies that for all borelian sets $E\subset \R^N$, one has:
\begin{equation}\label{eq99}
\nu(E)\LEQ c_{\mu}\mu(E)^{\frac{p^*}{p_{+}}}.
\end{equation}
Since the set $D=\{x\in \R^N:\mu(\{x\})>0\}$ is at most countable because $\mu\in\calM(\R^N)$, therefore $D=\{x_{j},\ j\in J\}$ and $b_{j}\dot=\mu(\{x_{j}\})$ satisfies $\DST\mu\GEQ\sum_{j\in J } b_{j} \delta_{x_{j}}$.\\
\ \\
Relation (\ref{eq99}) implies that $\nu$ is absolutely continuous with respect to $\mu$, {\it i.e.,} 
$\nu<\!\!<\mu$ and 
$$\frac{\nu\big(B(x,r)\big)}{\mu\big(B(x,r)\big)}\LEQ c_{\mu}\mu\big(B(x,r)\big)^{\frac{p^*}{p_{+}}-1},$$
provided that $\mu\big(B(x,r)\big)\neq0$ (remember that $p^*>p_{+}$). Thus, we have :
$$\nu(E)=\int_{E}\lim_{r\to0}\frac{\nu\big(B(x,r)\big)}{\mu\big(B(x,r)\big)}d\mu(x),$$
and
$$D_{\mu}\nu(x)=\lim_{r\to0}\frac{\nu\big(B(x,r)\big)}{\mu\big(B(x,r)\big)}=0,\ \mu\ a.e. \; 
\mbox{on} \; \R^N \setminus D.$$
Setting $a_{j}=D_{\mu}\nu(x_{j})b_{j}$, relation (\ref{eq99}) implies that $\nu$ has only atoms that are given by $\{x_{j}\}$, that we have already get.\\ 
\ \\
Let $\f\in C_{c}^\infty(\R^N),\ \f(x_{j})=1,\ \Norm\f_{\infty}=1$. Then, 
using statement 1. of this corollary and relation ({\ref{eq1000}), we have
\begin{equation}
Sa_{j}^{\frac{p_{+}}{p^*}}\LEQ
S\left(\int\norm\f^{p^*}d\nu
\right)^{\frac{p_{+}}{p^*}}
\LEQ\sum_{i=1}^N\frac1{p_{i}}\int\norm\f^{p_{i}}d\mu_{i}.
\end{equation}
We shall consider $\phi\in C_{c}^\infty(\R^N),\ 0\LEQ\phi\LEQ1$, 
support($\phi) \subset B(0,1)$, $\phi(0)=1$. We fix $j\in J$  and set 
$x_{j}=(x_{j,1},\ldots,x_{j,N}),\ q_{i}=\DST\frac{p_{i}p^*}{p^*-p_{i}},\ i=1,\ldots,N$.
Then $\DST\alpha_{i}\dot=\frac1{q_{i}}$ satisfy $\DST\sum_{k=1}^N\alpha_{k}-\alpha_{i}q_{i}=0$. For $\eps>0$, we define, for every $z\in\R^N,\ z=(z_{1},\ldots,z_{N})$:
\begin{equation}\label{eq989}
\phi_{\eps}(z)=\phi\left(\frac{z_{1}-x_{j,1}}{\eps^{\alpha_{1}}},\ldots,
\frac{z_{N}-x_{j,N}}{\eps^{\alpha_{N}}}\right).
\end{equation}
Thus we have :
\begin{equation}\label{1001}
\int\norm{\frac{\partial \phi_{\eps}}{\partial x_{i}}}^{q_{i}}
=\int\norm{\frac{\partial \phi}{\partial x_{i}}}^{q_{i}}(z)dz
\end{equation}
and then 
\begin{equation}\label{eq1002}
\int\norm  {\frac{\partial \phi_{\eps}}{\partial x_{i}}}^{p_{i}}
\norm v^{p_{i}}\LEQ
\left(\int\norm{\frac{\partial \phi}{\partial x_{i}}}^{q_{i}}dz\right)^{1-\frac{p_{i}}{p^*}}
\left(\int_{B(x_{j},\max_{i}\eps^{\frac1q_{i}})}
\norm v^{p^*}dz\right)^{\frac{p_{i}}{p^*}}
\xrightarrow[\eps\to0]{}0.
\end{equation} 

\begin{lem}\label{l12}
Let $x_{j}\in D$ and $\phi_{\eps}$ be the function defined above associated to $x_{j}$. Then :
$$Sa_{j}^{\frac{p_{+}}{p^*}}\LEQ\overline{\lim_{\eps\to0}} \, 
\overline{\lim_{n}}\sum_{i=1}^N\frac1{p_{i}}\int\phi_{\eps}^{p_{i}}\norm{\frac{\partial v_{n}}{\partial x_{i}}}^{p_{i}}dx.$$
\end{lem}
\PRF
Since $0\LEQ \phi_{\eps}\LEQ1$ then $\DST\int\phi_{\eps}^{p_{*}}\norm{v_{n}}^{p^*}dx\LEQ1$. From Corollary~\ref{c1l3} of Lemma~\ref{l3}, it follows
\begin{equation}\label{eq110}
S\left(\int\phi_{\eps}^{p_{*}}\norm{v_{n}}^{p^*}dx\right)^{\frac{p_{+}}{p^*}}
\LEQ\sum_{i=1}^N\frac1{p_{i}}\int\norm{\frac\partial{\partial x_{i}}(\phi_{\eps}v_{n})}^{p_{i}}.
\end{equation}
From relation (\ref{eq1002}), we have
\begin{equation}\label{eq111}
\lim_{\eps\to0}\int\norm{\frac{\partialÊ\phi_{\eps}}{\partial x_{i}}}^{p_{i}}\norm v^{p_{i}}dx=0.
\end{equation}
Since
\begin{equation}\label{eq112}
\lim_{n\to+\infty}\int\norm{\frac{\partialÊ\phi_{\eps}}{\partial x_{i}}}^{p_{i}}\norm {v_{n }-v}^{p_{i}}dx=0,
\end{equation}
then one has :
\begin{equation}\label{eq113}
\overline{\lim_{\eps\to0}}\,\overline{\lim_{n}}\sum_{i=1}^N\frac1{p_{i}}
\int\norm{\frac{\partialÊ}{\partial x_{i}}(\phi_{\eps}v_{n})}^{p_{i}}dx
=
\overline{\lim_{\eps\to0}}\, \overline{\lim_{n}}\sum_{i=1}^N\frac1{p_{i}}
\int\norm{\frac{\partial v_{n}Ê}{\partial x_{i}}}^{p_{i}}\phi_{\eps}^{p_{i}}dx
\end{equation}
From relations (\ref{eq110}) and (\ref{eq113}), knowing that 
$\norm{v_{n}}^{p^*}\rightharpoonup\norm v^{p^*}+\nu$ (see Lemma~\ref{l4}), we obtain
$$Sa_{j}^{\frac{p_{+}}{p^*}}
\LEQ\overline{\lim_{\eps\to 0}}\,\overline{\lim_{n}}
\sum_{i=1}^N\frac1{p_{i}}\int\phi_{\eps}^{p_{i}}\norm{\frac{\partial v_{n}}{\partial x_{i}}}^{p_{i}}dx.$$
\ \HF
\begin{lem}\label{l13}\ \\
Assume that $\DST\sum_{i=1}^N\frac1{p_{i}}\norm{\frac{\partial v_{n}}{\partial x_{i}}}^{p_{i}}\rightharpoonup \widetilde\mu$ in $\calM^+(\R^N)$. Then
\begin{enumerate}
\item For all $j\in J$, $\DST Sa_{j}^{\frac{p_{+}}{p^*}}\LEQ \lim_{\eps\to0}\widetilde\mu({\rm support}\phi_{\eps})$\\
 $($one has support $\phi_{\eps}\subset B(x_{j},\max_{i}\eps^{\frac1{q_{i}}}))$.
\item $\Norm{\widetilde \mu}\GEQ S\Norm \nu^{\frac{p_{+}}{p^*}}+P(\nabla v)$.
\item $S=\lim_{n \to +\infty} P(\nabla v_{n})\dot=\Norm{\widetilde\mu}+\mu_{\infty}\GEQ P(\nabla v)+S\Norm\nu^{\frac{p_{+}}{p^*}}+\mu_{\infty}$.
\end{enumerate}
\end{lem}
\PRF
From Lemma \ref{l12}, since $\phi_{\eps}^{p_{i}}\LEQ\phi_{\eps}$ and 
$$\overline{\lim_{n}}\sum_{i=1}^N\frac1{p_{i}}\int\phi_{\eps}^{p_{i}}\norm{\frac{\partial v_{n}}{\partial x_{i}}}^{p_{i}}dx\LEQ\int\phi_{\eps}d\widetilde \mu,$$
one obtains
\begin{equation}\label{eq114}
Sa_{j}^{\frac{p_{+}}{p^*}}\LEQ\overline{\lim_{\eps\to0}}\int\phi_{\eps}d\widetilde\mu
\LEQ\lim_{\eps\to0}\widetilde\mu\left(B(x_{j};\max_{1\LEQ i\LEQ N}\eps^{\frac1{q_{i}}})\right).
\end{equation}
This shows that $\{x_{j}\}_{j\in J}$ are all atomic points of $\widetilde \mu$ and since $\DST\sum_{i=1}^N\frac1{p_{i}}\norm{\frac{\partial v}{\partial x_{i}}}^{p_{i}}$ is orthogonal to the atomic part of $\widetilde \mu$, one deduces from relation (\ref{eq114}) that
\begin{equation}\label{eq115}
\widetilde \mu\GEQ S\sum_{j\in J}a_{j}^{\frac{p_{+}}{p^*}}\delta_{x_{j}}
+\sum_{i=1 }^N\frac1{p_{i}}\norm{ \frac{\partial v}{\partial x_{i}}}^{p_{i}}.
\end{equation}
This implies in particular that :
\begin{equation}\label{eq116}
\Norm{ \widetilde\mu}\GEQ 
S\sum_{j\in J}a_{j}^{\frac{p_{+}}{p^*}}+P(\nabla v).
\end{equation}
Since $\DST\frac{p_{+}}{p^*}<1$ one has
\begin{equation}\label{eq117}
\left(\sum_{j\in J}a_{j}\right)^{\frac{p_{+}}{p^*}}\LEQ
\sum_{j\in J}a_{j}^{\frac{p_{+}}{p^*}}.
\end{equation}
As $\DST \nu=\sum_{j\in J}a_{j}\delta_{x_{j}}$, it holds
\begin{equation}\label{eq118}
\Norm\nu=\sum_{j\in J}a_{j},
\end{equation}
which means, combining relations (\ref{eq116}) to (\ref{eq118}), that :
$$\Norm{\widetilde \mu}\GEQ S\Norm\nu^{\frac{p_{+}}{p^*}}+P(\nabla v).$$
For the last statement, we argue as before:
\begin{eqnarray*}
S&=&\lim_{n}P(\nabla v_{n})
\\&=&\lim_{R\to+\infty}\lim_{n}\int_{\R^N}(1-\psi_{R})
\sum_{i=1}^N\frac1{p_{i}}\norm{\frac{\partial v_{n}}{\partial x_{i}}}^{p_{i}}dx\\
&\ &+\lim_{R\to+\infty}\overline{\lim_{n}}\int\psi_{R}\sum_{i=1}^N
\frac1{p_{i}}\norm{\frac{\partial v_{n}}{\partial x_{i}}}^{p_{i}}dx,
\end{eqnarray*}
where $\psi_{R}=1$ on $\norm x>R+1$, $0\LEQ\psi_{R}\LEQ1$, 
$\psi_{R}=0$ if $\norm x<R$, $\psi_{R}\in C(\R)$.\\
By the definition of $\widetilde\mu$, one has :
$$\lim_{R\to+\infty}{\lim_{n}}\int(1-\psi_{R})\sum_{i=1}^N
\frac1{p_{i}}\norm{\frac{\partial v_{n}}{\partial x_{i}}}^{p_{i}}dx=
\lim_{R}\int(1-\psi_{R})d\widetilde\mu=\Norm{\widetilde\mu},$$
and (see Lemma \ref{l8}):
$$\lim_{R\to+\infty}\overline{\lim_{n}}\int\psi_{R}\sum_{i=1}^N
\frac1{p_{i}}\norm{\frac{\partial v_{n}}{\partial x_{i}}}^{p_{i}}dx=\mu_{\infty},$$
thus, by the preceding statements:
$$S=\Norm{\widetilde \mu}+\mu_{\infty}
\GEQ P(\nabla v)+S\Norm\nu^{\frac{p_{+}}{p^*}}+\mu_{\infty}.$$
\ \HF\\
\begin{lem}\label{l14}\ \\
If $\Norm v_{p^*}<1$ then $\Norm\nu=1,\ \nu_{\infty}=0\ and\ v=0.$
\end{lem}
\PRF
From Lemma \ref{l10}, we know that 
$$S\nu_{\infty}^{\frac{p_{+}}{p^*}}\LEQ\mu_{\infty}.$$
And by Corollary  \ref{c1l3} of Lemma \ref{l3}, we have
$$S\Norm v^{p_{+}}_{p^*}\LEQ P(\nabla v).$$
From the last statement of Lemma \ref{l13} and the above inequalities we deduce that :
$$S\GEQ S\Big((\Norm v^{p^*}_{p^*})^{\frac{p_{+}}{p^*}}+\Norm \nu^{\frac{p_{+}}{p^*}}
+\nu_{\infty}^{\frac{p_{+}}{p^*}}\Big).$$
Thus we obtain, due to Lemma \ref{l11}, that
$$\Big((\Norm v^{p^*}_{p^*})^{\frac{p_{+}}{p^*}}+\Norm \nu^{\frac{p_{+}}{p^*}}
+\nu_{\infty}^{\frac{p_{+}}{p^*}}\Big)\LEQ1
=\Big(\Norm v^{p^*}_{p^*}+\Norm \nu
+\nu_{\infty}
\Big)^{\frac{p_{+}}{p^*}}.$$
Using the inequality
$$\Big(\Norm v^{p^*}_{p^*}+\Norm \nu+\nu_{\infty}\Big)^{\frac{p_{+}}{p^*}}
\LEQ\Norm v_{p^*}^{\frac{p_{+}}{p^*}}+\Norm\nu^{\frac{p_{+}}{p^*}}+\nu_{\infty}^{\frac{p_{+}}{p^*}},$$
we get
$$\Norm v_{p^*}^{\frac{p_{+}}{p^*}}+\Norm\nu^{\frac{p_{+}}{p^*}}+\nu_{\infty}^{\frac{p_{+}}{p^*}}=\Big(\Norm v^{p^*}_{p^*}
+\Norm \nu
+\nu_{\infty}\Big)^{\frac{p_{+}}{p^*}}.$$
It follows that $\Norm v^{p^*}_{p^*},\ \Norm\nu$ and $\nu_{\infty}$ are equal either to 0 
or to 1. But using the fact that $\nu_{\infty}\DST\LEQ\frac12$, since $\DST\int_{B(0,1)}\norm{v_{n}}^{p^*}dx=\frac12,$ we conclude that $\nu_{\infty}=0$, 
$\Norm v_{p^*}<1$ (by our assumption) so that $v=0$ and thus $\Norm \nu=1$.\HF\\
\begin{lem}\label{l15}\ \\
If $\Norm v_{p^*}<1$ then the measure $\nu$ is concentrated at a single point 
$z=x_{i_{0}}$.
\end{lem}
\PRF
Since 
$$S=\Norm{\widetilde \mu}+\mu_{\infty}
\GEQ S\sum_{j\in J}a_{j}^{\frac{p_{+}}{p^*}},$$ (see relation(\ref{eq116})) and $1=\Norm \nu=\DST\sum_{j\in J}a_{j}$, we then have :
 $$\left(\sum_{j\in J}a_{j}\right)^{\frac{p_{+}}{p^*}}\GEQ
\sum_{j\in J}a_{j}^{\frac{p_{+}}{p^*}}
\GEQ \left(\sum_{j\in J} a_{j} \right)
^{\frac{p_{+}}{p^*}}.$$
Thus the $a_{j}$ are equal either to zero or to 1 that is, there is only one index 
$i_{0}$ such that $a_{i_{0}}=1$ and $a_{j}=0$ for 
$j\neq {i_0}:\nu=a_{i_{0}}\delta_{x_{i_{0}}}$.\HF\\
\  \\
{\bf End of the proof of the main Lemma :}\\
If $\Norm{v}_{p^*}<1$ thus $\nu$ concentrates at $x_{i_{0}}$ and $\Norm\nu=1$. 
On the other hand we have 
$\DST\frac12=\sup_{y\in\R^N}\int_{B(y,1)}\norm{v_{n}}^{p^*}
\GEQ\int_{B(x_{i_{0}},1)}
\norm{v_{n}}^{p^*}dx\to\Norm\nu=1$, which is impossible, we conclude then that $\Norm v_{p^*}=1$.\HF\\

Consequently, the function $v$ is a (non trivial) extremal function that can be chosen nonnegative (replacing $v$ by $\norm v$). \\ \ \\
{\bf End of the proof of Theorem \ref{t1} :}\\
From usual Lagrange multiplier rule, there is $\lambda_{0}>0$, such that :
$$-\sum_{i=1}^N\frac\partial{\partial x_{i}}\left(
\norm{\frac{\partial v}{\partial x_{i}}}^{p_{i}-2}\frac{\partial v}{\partial x_{i}}
\right)=\lambda_{0 }v^{p^*-1} \;\; \mbox{in} \;\; \calD^{1,\vp}(\R^N)'.$$

A similar rescaling  argument used above 
(say $v( \lambda_{0}  ^{-\frac1{p_{1}}}x_{1},\ldots,\lambda_{0} ^{-\frac 1{p_{N}}}x_{N})$~) 
gives the result.\HF\\

The multiplicity of solutions comes directly from Lemma~\ref{l2}, that is :
\begin{lem}\label{l14000}: \\
Let $\alpha\in\R,\ \alpha_{i}=\alpha\DST\frac{p^*}{p_{i}}-\alpha,\ i=1,\ldots,N$ and $u\in\calS$. Then, for all $\lambda\in\R_{+}^*$ for all $z=(z_{1},\ldots,z_{N})\in\R^N$, the function defined by
$$u^{\lambda,z}(x)=\lambda^\alpha u(\lambda^{\alpha_{1}}x_{1}+z_{1},\ldots,
\lambda^{\alpha_{N}}x_{N}+z_{N}),$$
with $x=(x_{1},\ldots,x_{N})$ belongs to $\calS$.
\end{lem}
\PRF
It is the same as for Lemma \ref{l2} using a direct computation.\\

\section{Some properties of the solutions of (1)}
We want to show first the :
\begin{prop}\label{p1}\ \\
Any nonnegative solution $u$ being in $\calD^{1,\vp}(\R^N)$ of (1) belongs to $L^q(\R^N)$ for all $p^*\LEQ q<+\infty$.
\end{prop}
\PRF
We follow the proof of \cite {FGK}. Let $a>0$.
Let $j$ be fixed in $\{1,\ldots,N\}$, for $L>0$ (large) we define $\f_{j,L}\dot=u\min[u^{ap_{j}},L^{p_{j}}]\in\calD^{1,\vp}(\R^N)$ and for all $i$
\begin{equation}\label{eq8000}
\norm{\partial _{i}u}^{p_{j}-2}\partial_{i}u\partial_{i}\f_{j,L}\GEQ\min[u^{ap_{j}},L^{p_{j}}]\norm{\partial _{i}u}^{p_{j}}\quad a.e,
\end{equation}
and {\begin{equation}\label{eq8001}
\norm{\partial_{i}(u\cdot\min[u^{a},L])}^{p_{j}}\LEQ (a+1)^{p_{j}}\min[u^{ap_{j}},L^{p_{j}}]\norm{\partial_{i}u}^{p_{j}} \quad a.e.
\end{equation}
Choosing $\f_{j,L}$ as a test function, one has :
\begin{eqnarray}\label{eq8002}
\int_{\R^N}\min[u^{ap_{j}},L^{p_{j}}]\norm{\partial_{j}u}^{p_{j}}dx
&\LEQ&\sum_{i=1}^N\int_{\R^N}\norm{\partial_{i}u}^{p_{i}-2}\partial_{i}u\partial_{i}\f_{j,L}dx
\nonumber\\
&=&\int_{\R^N}u^{p^*}\min[u^{ap_{j}},L^{p_{j}}]dx.
\end{eqnarray}
Introducing $k>0$, one has :
\begin{equation}\label{eq8003}
\int_{\R^N}{u^{p^*}}\min[u^{ap_{j}},L^{p_{j}}]dx\LEQ k^{ap_{j}}\int_{\R^N}u^{p^*}dx
+\int_{u\GEQ k}u^{p^*}\min[u^{ap_{j}},L^{p_{j}}]dx.
\end{equation}
Writing that :
\begin{equation}\label{eq8004}
\int_{u\GEQ k}u^{p^*}
\min[u^{ap_{j}},L^{p_{j}}]dx=\int_{u\GEQ k}u^{p^*-p_{j}}u^{p_{j}}
\left(\min[u^{a},L]\right)^{p_{j}}dx.
\end{equation}
The H\"older inequality applied to the right hand side of relation ({\ref{eq8004}}) shows that :
\begin{equation}\label{eq8005}
\int_{u\GEQ k}u^{p^*}\min[u^{ap_{j}},L^{p_{j}}]dx
\LEQ\left(\int_{u\GEQ k}u^{p^*}dx\right)^{1-\frac{p_{j}}{p^*}}
\left(\int_{\R^N}\left( u \min[u^{a},L]\right)^{p^*}\right)^{\frac{p_{j}}{p^*}}.
\end{equation}
By the Troisi's inequality (see Lemma \ref{l1})
\begin{equation}\label{eq8006}
\left(\int_{\R^N}
\left(u\min[u^{a },L] \right)^{p^*}
\right)^{\frac1{p^*}}
\LEQ c\sum_{i=1}^N
\left(\int_{\R^N}
 \norm{\partial_{i}
 \left(u\min[u^{a},L]\right)}^{p_{i}}
 \right)^{\frac{1}{p_{i}}}
\end{equation}
Setting $I_{i}
=\DST
\left(\int\norm{\partial_{i}(u\min[u^{a},L])}^{p_{i}}\right)^{\frac1{p_{i}}}$, 
$\eps_{k}=\DST\int_{u\GEQ k}u^{p^*}dx$, 
 relations (\ref{eq8001}) to (\ref{eq8006}), lead to :
\begin{eqnarray*}
\int\norm{\partial_{j}(u\cdot\min[u^{a},L])}^{p_{j}}dx
&\LEQ&(a+1)^{p_{j}}\int\min[u^{ap_{j}},L^{p_{j}}]\norm{\partial_{j}u}^{p_{j}}dx\\
&\LEQ&(a+1)^{p_{j}}k^{ap_{j}}
\left(\int u^{p^*}dx\right)\\&&
\!\!\!\!\!+c(a+1)^{p_{j}}
\eps_{k}^{1-\frac{p_{j}}{p^*}}   
\left[\sum_{i=1}^N\left(\int\norm{\partial_{i}(u\min[u^{a},L])}^{p_{i}}\right)^{\frac1{p_{i}}}\right]^{p_{j}}.
\end{eqnarray*}
Thus, for all $j$ :
\begin{equation}\label{eq8007}
I_{j}\LEQ(a+1)k^{a}
\left(\int u^{p^*}dx\right)^{\frac1{p_{j}}}+c(a+1)\eps_{k}
^{\frac1{p_{j}}-\frac1{p^*}}\left(\sum_{i=1}^NI_{i}\right)
\end{equation}

The relation(\ref{eq8007}) infers :
\begin{equation}\label{eq8008}
\sum_{j=1}^N I_{j}\LEQ(a+1)k^{a}\left(\sum_{j=1}^N\Norm u_{p^*}^{\frac{p^*}{p_{j}}}\right)
+c(a+1)\left(\sum_{j=1}^N\eps_{k}^{\frac1{p_{j}}-\frac1{p^*}}\right)
\left(\sum_{i=1}^NI_{i}\right).
\end{equation}
Since $\DST\lim_{k\to+\infty}\sum_{j=1}^N\eps_{k}^{\frac1{p_{j}}-\frac1{p^*}}=0$, 
there exists $k_{a}>0$ such that for all $k\GEQ k_{a}$, such that 
$\DST c(a+1)\sum_{j=1}^N\eps_{k}^{\frac1{p_{j}}-\frac1{p^*}}
\LEQ\frac12$. Thus relation (\ref{eq8008}) infers then 
$$\DST\sum_{i=1}^NI_{j}\LEQ 2(a+1)k^{a}\sum_{j=1}^N\Norm u_{p^*}^{\frac{p^*}{p_{j}}}, \ for\ k\GEQ k_{a}.$$
By the Troisi's inequality, one has :
$$\Norm{ u\cdot \min[u^a,L]}_{L^{p^*}}\LEQ c\sum_{j=1}^NI_{j}
\LEQ2c(a+1)k^{a}\sum_{j=1}^N\Norm u_{p^*}^{\frac{p^*}{p_{j}}}.$$
Letting $ L\to+\infty$, one has :
$$\DST\Norm{u^{a+1}}_{L^{p^*}}\LEQ 2c(a+1)k^{a}\sum_{j=1}^N 
\Norm u_{p^*}^{\frac{p^*}{p_{j}}}.$$
Let $q=(a+1){p^*}$, then we obtain the result.\HF

\begin{prop}\label{p2}
Any nonnegative solution $u$ being in 
$\calD^{1,\vp}(\R^N)$ of (1) belongs to $L^\infty(\R^N)$.
 Moreover, there exists a number $\tau_{0}$ depending only on $p_{j}$, $N$ such that $$\Norm u_{p^*}\GEQ\tau_{0}>0,\ for \ u \hbox{ non trivial}.
 $$  
\end{prop}
\PRF
For $u\GEQ0$ solution of (1), we set $A_\tau=\{x\in\R^N,\ u(x)\GEQ \tau\}$ 
and $\norm{A_{\tau}}$ its Lebesgue measure. 
Since $p^*>p_{+}$, one can choose $q>p^*$ so that 
$$\eps\dot=-\frac1{p^*}+\left(1-\frac{p^*}q\right)\left(1-\frac1{p^*}\right)\frac1{p_{+}-1}>0.$$
Let $\f_k=(u-k)_{+}$, for $k>0$ fixed. Chosing this function as a test function and using proposition \ref{p1}, one has :
\begin{equation}\label{eq8009}
\sum_{i=1}^N\Norm{\frac{\partial \f_{k}}{\partial x_{i}}}_{p_{i}}^{p_{i}}
=\int u^{p^*-1}(u-k)_{+}\LEQ c_{1}\norm{A_{k}}^{\left(1-\frac{p^*}q\right)\left(1-\frac1{p^*}\right)}\Norm{\f_{k}}_{p^*},
\end{equation} 
with $c_{1}=\Norm u_{q}^{p^*-1}.$\\
Since $\DST\Norm{\f_{k}}_{p^*}\LEQ\Norm u_{p^*}$, 
thus the corollary~\ref{c1l3} of Lemma~\ref{l3} and relation (\ref{eq8009}) imply~:
\begin{equation}\label{eq8010}
\Norm{\f_{k}}^{p_{+}}_{p^*}\LEQ c_{2} \sum_{i=1}^N
\Norm{\frac{\partial \f_{k}}{\partial x_{i}}}^{p_{i}}_{p_{i}}
\LEQ c_{3}\norm{A_{k}}^{\left(1-\frac{p^*}q\right)\left(1-\frac1{p^*}\right)}
\Norm{\f_{k}}_{p^*},
\end{equation}
with $c_{2}=\DST\frac1{S\cdot p_{-}}\Max_{1\LEQ j\LEQ N}\left(\Norm u_{p^*}^{p_{+}-p_{j}}\right),\ 
c_{3}=c_{1}c_{2}.$\\
Thus, 
\begin{equation}\label{eq8011}
\Norm{\f_{k}}_{p^*}\LEQ c_{4}\norm{A_{k}}^{\frac1{p_{+}-1}\left(1-\frac{p^*}q\right)\left(1-\frac1{p^*}\right)}.
\end{equation}
with $c_{4}=\DST c_{3}^{\frac1{p_{+}-1}}$.
By Cavalieri's principle, H\"older inequality and relation(\ref{eq8011}), one has,  for all $k>  {0}$:
\begin{equation}\label{eq8012}
\int_{k}^{+\infty}\norm{A_{\tau}}d\tau
=\int_{\R^N}(u-k)_{+}(x)dx
\LEQ\norm{A_{k}}^{1-\frac1{p^*}}\Norm{\f_{k}}_{p^*}
\LEQ c_{4}\norm{A_{k}}^{1+\eps}.
\end{equation}
This last relation is a Gronwall inequality, which shows that $\forall k>0$
\begin{equation}\label{eq38}
\Norm u_{\infty}\LEQ k+\frac{1+\eps}\eps\Norm{(u-k)_{+}}_{1}^{\frac\eps{1+\eps}}
c_{4}^{\frac1{1+\eps}}.
\end{equation}
Setting
$$\gamma=(p^*-1)\frac\eps{1+\eps},\ b_{0}=\frac{1+\eps}\eps\Norm u_{p^*}^{\frac{\eps p^*}{1+\eps}}c_{4}^{\frac1{1+\eps}},$$
and noticing that 
$$\Norm{(u-k)_{+}}_{1}\LEQ\frac{\Norm u_{p^*}^{p^*}}{k^{p^*-1}},$$
thus relation(\ref{eq38}) becomes :
\begin{equation}\label{eq39}
\Norm u_{\infty}\LEQ \Inf_{k>0}\left[k+\frac{b_{0}}{k^\gamma}\right]
=(\gamma+1)\gamma^{-\frac\gamma{\gamma+1}}b_{0}^{\frac1{1+\gamma}}.
\end{equation}
Separating the contribution of $\Norm u_{q}$ and $\Norm u_{p^*}$, we have a continuous  map $\Lambda:\R_{+}\to\R_{+}$ and constants $c_{5}>0$ and $\beta$ depending only on $p_{+},\ p_{*}$ so that
\begin{equation}\label{eq40}
\Norm u_{\infty}\LEQ c_{5}\Norm u_{q}^\beta\Lambda(\Norm u_{p^*}),
\end{equation}
with $\DST\beta=\frac{p^*-1}{(p_{+}-1)(1+\eps)(1+\gamma)}$,\  
$\DST \Lambda(\sigma)=\left[\sigma^{\eps p^*}
\Max_{1\LEQ j\LEQ N}(\sigma^{p_{+}-p_{j}})
\right] 
^{\frac1{(1+\eps)(1+\gamma)}}$. \\
 Thus, from relation (\ref{eq40}),we deduce 
\begin{equation}\label{eq41}
\Norm u_{\infty}^{1-\beta(1-\frac{p^*}q)}\LEQ c_{5}\Norm u_{p^*}^{\beta\frac{p^*}q}
\Lambda(\Norm u_{p^*})\ for \ u\not\equiv0.
\end{equation}
But the number $\DST\kappa\dot=1-\beta\left(1-\frac{p^*}q\right)=0$, so relation (\ref{eq41}) implies that there is a number $\tau_{0}>0$ depending only $p_{j},\ p^*$ such that $\Norm u_{p^*}\GEQ\tau_{0}>0$.  \HF\\
\ \\
{\bf Acknowledgment :} The authors would like to thank the referee for his/her valuable comments.

%
%
{

}
\end{document}